\documentclass{amsart}

\usepackage{eucal}
\input{xypic}
\xyoption{all}

\begin{document}

\newcommand{\C}{{\mathbb C}}
\newcommand{\Diff}{{\rm Diff}}
\newcommand{\e}{{\bf e}}
\newcommand{\Emb}{{\rm Emb}}
\newcommand{\gr}{{\rm Grass}}
\newcommand{\G}{{\mathbb G}}
\newcommand{\half}{{\textstyle{\frac{1}{2}}}}
\newcommand{\aitch}{{\mathbb H}}
\newcommand{\kO}{{\rm kO}}
\newcommand{\M}{{\overline{\mathcal M}}}
\newcommand{\MU}{{\rm MU}}
\newcommand{\MSpin}{{\rm MSpin}}
\newcommand{\Q}{{\mathbb Q}}
\newcommand{\R}{{\mathbb R}}
\newcommand{\res}{{\rm res}}
\newcommand{\SU}{{\rm SU}}
\newcommand{\Spin}{{\rm Spin}}
\newcommand{\tr}{{\rm Trace}}
\newcommand{\T}{{\mathbb T}}
\newcommand{\thom}{{\rm Th}}
\newcommand{\x}{{\surd x}}
\newcommand{\Z}{{\mathbb Z}}

\title {Heisenberg groups and algebraic topology} 
\author{Jack Morava}
\address{Department of Mathematics, Johns Hopkins University, Baltimore,
Maryland 21218}
\email{jack@math.jhu.edu}
\thanks{The author was supported in part by the NSF}
\date {6 January 2003}
\begin{abstract} 
We study the Madsen-Tillmann spectrum $\C P^\infty_{-1}$ as a quotient 
of the Mahowald pro-object $\C P^{\infty}_{-\infty}$, which is closely 
related to the Tate cohomology of circle actions. That theory has an
associated symplectic structure, whose symmetries define the Virasoro 
operations on the cohomology of moduli space constructed by Kontsevich 
and Witten. \end{abstract}

\maketitle
\section {Introduction} \bigskip

\noindent
A sphere $S^n$ maps essentially to a sphere $S^k$ only if $n \geq k$, 
and since we usually think of spaces as constructed by attaching cells, 
it follows that algebraic topology is in some natural sense upper-triangular, 
and thus not very self-dual: as in the category of modules over the
mod $p$ group ring of a $p$-group, its objects are built by iterated 
extensions from a small list of simple ones. \bigskip

\noindent
Representation theorists find semi-simple categories more
congenial, and for related reasons, physicists are happiest in
Hilbert space. This paper is concerned with some remarkable
properties of the cohomology of the moduli space of Riemann
surfaces discovered by physicists studying two-dimensional
topological gravity (an enormous elaboration of conformal field
theory), which appear at first sight quite unfamiliar. Our
argument is that these new phenomena are forced by the physicists' 
interest in self-dual constructions, which leads to objects which are
(from the point of view of classical algebraic topology) very
large [1 \S 2]. \bigskip

\noindent
Fortunately, equivariant homotopy theory provides us with tools
to manage these constructions. The first section below is a
geometric introduction to the Tate cohomology of the circle group; 
the conclusion is that it possesses an intrinsic symplectic module 
structure, which pairs positive and negative dimensions in a way 
very useful for applications. Section two studies operations on this 
(not quite cohomology) functor, and exhibits the action of an algebraic 
analog of the Virasoro group on it. The third section relates rational 
Tate cohomology of the circle to that of the infinite loopspace
$Q\C P^\infty_+$ considered by Madsen and Tillmann in recent 
work on Mumford's conjecture. \bigskip

\noindent
I owe thanks to many people for help with the ideas in this paper, but
it is essentially a collage of a lifetime's conversations with 
Graeme Segal, who more or less adopted me when we were both very young.

\section {Geometric Tate cohomology}

\noindent
{\bf 2.1} Let $G$ be a compact Lie group of dimension $d$. We
will be concerned with a cobordism category of smooth compact 
$G$-manifolds, with the action free on the boundary: this can 
be regarded as a categorical cofiber for the forgetful functor from 
manifolds with free $G$-action to manifolds with unrestricted action. 
Under reasonable assumptions this cofiber category is closed under
Cartesian products (given the diagonal action). \bigskip

\noindent
If $E$ is a geometric cycle theory (eg stable homotopy, or classical 
homology) then the graded $E$-bordism group of free $G$-manifolds is 
isomorphic to $E_{*+d}(BG_+)$. On the other hand, the homotopy quotient of a
$G$-manifold is a bundle of manifolds over the classifying space
$BG$, and Quillen's conventions [22] associate to such a
thing, a class in the graded cobordism group $E^{-*}(BG_+)$. The
forgetful functor from free to unrestricted $G$-manifolds defines
a long exact sequence 
\[
\dots \to E_{*-d}(BG_+) \to E^{-*}(BG_+) \to t^{-*}_G(E) \to 
E_{*-d-1}(BG_+) \to \dots
\]
which interprets the relative groups as the (Tate-Swan [28]) $E$-bordism of
manifolds with $G$-action free on the boundary. The geometric
boundary homomorphism
\[
\partial_E : t^{-*}_G(E) \to E_{*-d-1}(BG_+) \to E_{*-d-1}
\]
sends a manifold with $G$-free boundary to the quotient of that
action on the boundary; it will be useful later. \bigskip
 
\noindent
Remarks: \medskip

\noindent
1) $t_G(E)$ is a ring-spectrum if $E$ is; in fact it is an $E$-algebra, 
at least in some naive sense; \medskip

\noindent
2) tom Dieck stabilization [4] extends this geometric bordism
theory to an equivariant theory; \medskip

\noindent
3) the functor $t_G$ sends cofibrations to cofibrations, but it
lacks good limit properties: it is defined by a kind of hybrid
of homology and cohomology, and Milnor's limit fails. In
more modern terms [11], the construction sends a $G$-spectrum $E$ to
the equivariant function spectrum $[EG_+,E] \wedge \tilde EG_+$.
\medskip

\noindent
4) The eventual focus of \S 2 is the case when $G$ is the circle 
group $\T$, and $E$ is ordinary cohomology: this is closely related 
to cyclic cohomology [2], but I don't know enough about that subject to 
say anything useful. \bigskip

\noindent
{\bf 2.2} Suppose now that $E$ is a general complex-oriented
ring-spectrum; then $E^*(B\T_+)$ is a formal power series ring
generated by the Euler (or first Chern) class $\e$. If $E_*$
is concentrated in even degrees, then the cofiber sequence above
reduces for dimensional reasons to a short exact sequence
\[
0 \to E^*(B\T_+) = E^*[[\e]] \to t^*_\T E = E^*((\e)) \to 
E_{-*-2}(B\T_+) \to 0 
\]
with middle group the ring of formal Laurent series in $\e$. By
a lemma  of [12 \S 2.4] we can think of $t^*_\T E$ as the
homotopy groups of a pro-object
\[
S^2\C P^\infty_{-\infty} \wedge E := \{ S^2\thom(-k\eta) \wedge E
\;|\; k > -\infty \}
\]
in the category of spectra, constructed from the filtered vector
bundle
\[
\dots (k-1) \eta \subset k \eta \subset (k+1) \eta \subset \dots 
\]
defined by sums of copies of the tautological line bundle over 
$\C P^\infty \cong B\T$, as discussed in the appendix to [6] (see
also [18]). More precisely: the Thom spectrum can be taken to be 
\[
\thom (-k \eta) := \lim_n  S^{-2(n+1)k} \C P_n^{k \eta^\perp} \;,
\]
where $\eta^\perp$ is the orthogonal complement to the canonical
line $\eta$ in $\C^{n+1}$.] \bigskip

\noindent
When $E$ is {\bf not} complex-orientable, $t_\T E$ can behave
very differently: the Segal conjecture for Lie groups implies
that, up to a profinite completion [10], 
\[
t_\T S^0 \sim S^0 \vee S[\prod B\T/C]
\]
where the product runs through proper subgroups $C$ of $\T$ (and
$S$ denotes suspension). \bigskip

\noindent
In the universal complex-oriented case, the class $\e^{-1} \in t_\T^{-2} \MU$ 
is represented geometrically by the unit disk in $\C$
with the standard action of $\T$ as unit complex numbers; more
generally, the unit ball in $\C^k$ represents $\e^{-k}$. The
geometric boundary homomorphism sends that $\T$-manifold to $\C
P_{k-1}$; this observation can be restated, using Mishchenko's
logarithm, as the formula
\[
\partial_E (f) = \res_{\e = 0} \; f(\e) \; d \log_{\MU}(\e) : t^*_\T
MU \to \MU_{-*-2} \;,
\]
where the algebraic residue homomorphism 
\[
\res_{\e = 0} : \MU^{*-2}((\e)) \to \MU^*
\]
is defined by $\res_{\e=0}\; \e^k \; d\e = \delta_{k+1,0}$, cf.
[19, 21, 29]. \bigskip

\noindent
{\bf 2.3} The relative theory of manifolds with free group action
on the boundary alone defines bordism groups $\tau^*_G (E)$
analogous to a truncation of Tate cohomology, with useful
geometric applications. In place of the long exact sequence
above, we have 
\[
\dots \to E^{-*}(S^0) \to \tau^{-*}_G(E) \to E_{*-d-1}(BG_+) \to
\dots 
\]
compatible with a natural transformation $t^*_G(E) \to
\tau^*_G(E)$ which forgets the interior $G$-action. In our case
(when $E$ is complex-oriented), this is just the $E$-homology of 
the collapse map
\[
\C P^\infty_{-\infty} \to \C P^\infty_{-1} := \thom(-\eta)
\]
defined by the pro-spectrum in the previous paragraph. \bigskip

\noindent
A Riemann surface with geodesic boundary is in a natural way an
orientable manifold with a free $\T$-action on its boundary, and
a family of such things, parametrized by a space $X$, defines an
element of 
\[
\tau^{-2}_\T \MU(X_+) \cong [X, \C P^\infty_{-1} \wedge \MU] \;.
\]
The Hurewicz image of this element in ordinary cohomology is the
homomorphism 
\[
H^*(\C P^\infty_{-1},\Z) \to H^*(X_+,\Z)
\]
defined by the classifying map of Madsen and Tillmann, which will
be considered in more detail below. \bigskip 
 
\section{Automorphisms of classical Tate cohomology}

\noindent
{\bf 3.1} There are profound analogies -- and differences --
among the Tate cohomology rings of the groups $\Z/2\Z, \T$, and
$\SU(2)$ [3]. A property unique to the circle is the 
existence of the nontrivial involution $I : z \to z^{-1}$. \bigskip

\noindent
When $E$ is complex oriented, the symmetric bilinear form
\[
f,g \mapsto (f,g) = \partial_E (fg)
\]
on the Laurent series ring $t^*_\T E$ is nondegenerate, and the
involution on $\T$ defines a symplectic form
\[
\{ f,g \} = (I(f),g)
\]
which restricts to zero on the subspace of elements of positive
(or negative) degree. This Tate cohomology thus has an intrinsic
inner product, with canonical polarization and involution.
\bigskip 

\noindent
This bilinear form extends to a generalized Kronecker pairing 
\[
t_\T E^*(X) \otimes_E t_\T E_*(X) \to E_{*-2}
\]
which can be interpreted as a kind of Spanier-Whitehead duality
between $t_\T E_*$, viewed as a pro-object as in \S 2.1, and the
direct system $\{E^*(\thom(-k \eta)) \;|\; k > -\infty \}$ defined
by the cohomology of that system. This colimit again defines a Laurent
series ring, but this object is not quite its own dual: a shift of 
degree two intervenes, and it is most natural to think of the 
(non-existent) functional dual of $t_\T E$ as $S^{-2} t_\T E$. 
The residue map $t_\T E \to S^2 E$ can thus be understood as dual 
to the unit ring-morphism $E \to t_\T E$. \bigskip

\noindent
{\bf 3.2} The Tate construction is too large to be conveniently
represented, so the usual Hopf-algebraic approach to the study of
its automorphisms is technically difficult. Fortunately, methods
from the theory of Tannakian categories can be applied: we consider
automorphisms of $t_\T E$ as $E$ varies, and approximate
the resulting group-valued functor by representable ones. There
is no difficulty in carrying this out for a general complex-oriented 
theory $E$, but the result is a straightforward extension of the 
case of ordinary cohomology. \bigskip

\noindent
To start, it is clear that the group(scheme, representing the
functor
\[
A \mapsto \G_0(A) \; = \; \{ g(x) = \sum_{k \geq 0} g_k x^{k+1} \in
A[[x]] \;|\; g_0 = 1 \}
\]
on commutative rings $A$) of automorphisms of the formal line
acts as multiplicative natural transformations of the 
cohomology-theory-valued functor $A \mapsto t^*_\T HA$, with $g 
\in \G_0(A)$ sending the Euler class $\e$ to $g(\e)$. [I am treating 
these theories as graded by $\Z/2\Z$, with $A$ concentrated in degree
zero; but one can be more careful.] \bigskip

\noindent
Clearly $\G_0$ is represented by a polynomial Hopf algebra on
generators $g_k$, with diagonal
\[
(\Delta g)(x) = (g \otimes 1)((1 \otimes g)(x)) \;.
\]
However, $\G_0$ is a subgroup of a larger system $\G$ of natural
automorphisms, which is a colimit of representable functors
(though not itself representable): following [16], let 
\[
A \mapsto \G (A) \;  = \; \{ g(x) = \sum_{k \gg -\infty} g_k x^{k+1} \in
A((x)) \;|\; g_0 \in A^\times, g_k \in \surd A \; {\rm if} \; 0 >
k \}
\]
be the group of invertible {\bf nil-Laurent} series, ie Laurent
series with $g_0$ a unit, and $g_k$ nilpotent for negative $k$.
It is clear that $\G$ is a monoid, but in fact [20] it possesses
inverses. The Lie algebra of $\G$ is spanned by the derivations 
$x^{k+1} \partial_x, \; k \in \Z$: it is the algebra of vector 
fields on the circle. \bigskip

\noindent
{\bf 3.3} A related group-valued functor preserves the symplectic
structure defined above: to describe it, I will specialize even
further, and work over a field in whch two is invertible: $\R$,
for convenience. Thus let $\check \G$ be the (ind-pro)-algebraic
groupscheme defined by invertible nil-Laurent series over the
field $\R((\x))$ obtained from $\R((x))$ by adjoining a formal
square root of $x$, and let $\check \G_{\rm odd}$ denote the 
subgroup of {\bf odd} invertible series $\check g(\x) = - \check
g(-\x)$. The homomorphism 
\[
\check g \mapsto g(x) := \check g(\x)^2 : \check \G_{\rm odd} \to
\G
\]
is then a kind of double cover. \bigskip

\noindent
The functor $\check \G$ acts by symplectic automorphisms of the
module $\R((\x))$, given the bilinear form
\[
\langle u,v \rangle := \pi \; \res_{x=0} \; u(x) \; dv(x) 
\]
[27]; it is in fact a group of restricted symplectic
automorphisms of this module. The Galois group of
$\R((\x))/\R((x))$ defines a $\Z/2\Z$ - action, and the subgroup 
$\check \G_{\rm odd}$ preserves the subspace $\R((\x))_{odd}$
of odd power series. \bigskip

\noindent
{\bf Proposition} The linear transformation
\[
t^*_\T H\R \to \R((\x))_{odd}
\]
defined on normalized basis elements by
\[
\e^k \mapsto \gamma_{-k-\half}(x)
\]
(where $\gamma_s(x) = \Gamma(1+s)^{-1} x^s$ denotes a divided power),
is a dense symplectic embedding. \bigskip

\noindent   
Proof: We have 
\[
\{\e^k,\e^l \} = (-1)^k \; \res_{\e=0} \; \e^{k+l} \; d\e = (-1)^k
\delta_{k+l+1,0}
\]
while 
\[
\langle \gamma_s,\gamma_t \rangle = \res_{x=0} \; \gamma_s(x)
\gamma_{t-1}(x) \; dx = \frac{\pi}{\Gamma(t) \Gamma(1+s)}
\delta_{s+t,0} \;.
\]
The assertion then follows from the duplication formula for the Gamma
function. \bigskip

\noindent
The half-integral divided powers lie in $\Q((\x))$, aside from
distracting powers of $\pi$. The remaining rational coefficients
involve the characteristic `odd' factorials of 2D topological
gravity [8, 15], eg when $k$ is positive,
\[
\Gamma(k + \half) = (2k-1)!! \; 2^{-k} \surd \pi \;.
\]

\section {Symmetries of the stable cohomology of the Riemann moduli space}

\noindent
The preceding sections describe the construction of a polarized
symplectic structure on the Tate cohomology of the circle group.
The algebra of symmetric functions on the Lagrangian submodule 
\[
H^*(\C P^\infty_+,\Q) \subset t_\T H\Q
\]
of that cohomology can be identified with the homology of the
infinite loopspace 
\[
Q\C P^\infty_+ = \lim_n \Omega^n S^n \; \C P^\infty_+ \;;
\]
on the other hand, this module of functions admits a canonical 
action of the Heisenberg group associated to its defining symplectic
module [24 \S 9.5]. \bigskip

\noindent
The point of this paper is that the homology of this infinite
loopspace, considered in this way as a Fock representation,
manifests the Virasoro representation constructed by Witten and
Kontsevich on the stable cohomology of the moduli space of
Riemann surfaces, identified with $H_*(Q\C P^\infty_+,\Q)$
through work of Madsen, Tillman, and Weiss. Some of those results 
are summarized in the next two subsections; a more thorough 
account can be found in Michael Weiss's survey in these Proceedings.
The third section below discusses their connection with representation 
theory. \bigskip

\noindent
{\bf 4.1} Here is a very condensed account of one component of [17]: 
if $F \subset \R^n$ is a closed connected two-manifold embedded smoothly 
in a high-dimensional Euclidean space, its Pontrjagin-Thom construction 
$\R^n_+ \to F^{\nu}$ maps compactified Euclidean space to the Thom 
space of the normal bundle of the embedding. The tangent plane to $F$ 
is classified by a map $\tau : F \to \gr_{2,n}$ to the Grassmannian of 
oriented two-planes in $\R^n$, and the canonical two-plane bundle $\eta$ 
over this space has a complementary $(n-2)$-plane bundle, which I will 
call $(n - \eta)$. The normal bundle $\nu$ is the pullback along $\tau$ 
of $(n - \eta)$; composing the map induced on Thom spaces with the 
collapse defines the map
\[
\R^n_+ \to F^\nu \to \gr_{2,n}^{(n - \eta)} \;.
\]
The space $\Emb(F)$ of embeddings of $F$ in $\R^n$ becomes 
highly connected as $n$ increases, and the group $\Diff(F)$ of
orientation-preserving diffeomorphisms of $F$ acts freely
on it, defining a compatible family
\[
\R^n_+ \wedge_{\Diff} \Emb(F) \to \gr_{2,n}^{(n - \eta)}
\] 
which can be interpreted as a morphism
\[
B\Diff(F) \to \lim \; \Omega^n \gr_{2,n}^{(n - \eta)} \; := \; 
\Omega^\infty \C P^\infty_{-1}\;.
\]
This construction factors through a map
\[
\coprod_{g \geq 0} B\Diff (F_g) \to \Z \times B\Gamma_\infty^+ \to
\Omega^\infty \C P^\infty_{-1}
\]
of infinite loopspaces. Collapsing the bottom cell defines a
cofibration
\[
S^{-2} \to \C P^\infty_{-1} \to \C P^\infty_+
\]
of spectra; the fiber of the corresponding map
\[
\Omega^2 QS^0 \to \Omega^\infty \C P^\infty_{-1} \to Q\C
P^\infty_+
\]
of spaces has torsion homology, and the resulting composition
\[
\Z \times B\Gamma^+_\infty \to Q\C P^\infty_+ \sim QS^0 \times
Q\C P^\infty
\]
is a rational homology isomorphism which identifies Mumford's
polynomial algebra on classes $\kappa_i, \; i \geq 1$, with the
symmetric algebra on positive powers of $\e$. The rational
cohomology of $QS^0$ adds a copy of the group ring of $\Z$, which
can be interpreted as a ring of Laurent series in a zeroth
Mumford class $\kappa_0$. \bigskip

\noindent
The standard convention is to write $b_k$ for the generators 
of $H_*\C P^\infty_+$ dual to $\e^k$, and to use the same symbols for 
their images in the symmetric algebra $H_*(Q\C P^\infty_+,\Q)$. The Thom 
construction defines a map
\[
\C P^\infty \to \MU
\]
which extends to a ring isomorphism
\[
H_*(Q\C P^\infty,\Q) \to H_*(\MU,\Q) \;;
\]
sending the $b_k$ to classes usually denoted $t_k$, with $k \geq 1$;
but it is convenient to extend this to allow $k=0$. \bigskip

\noindent
{\bf 4.2} The homomorphism
\[
\lim \MU^{*+n-2}(\thom(n - \eta)) \to \MU^{*-2}(B\Diff(F)) 
\]
defined on cobordism by the Madsen-Tillmann construction sends the Thom 
class to a kind of Euler class: according to Quillen, the Thom class 
of $n - \eta$ is its zero-section, regarded as a map between manifolds. 
Its image is the class defined by the fiber product
\[
\xymatrix{
{Z_n} \ar[r] \ar[d]& {\gr_{2,n}} \ar[d]  \\
{\R^n_+ \wedge_{\Diff} \Emb(F)} \ar[r]& {\gr_{2,n}^{(n - \eta)}} \;;}
\]
this is the space of equivalence classes, under the action of
$\Diff(F)$, of pairs $(x,\phi)$, with $x \in \phi(F) \subset \R^n$
a point of the surface (ie, in the zero-section of $\nu$), and $\phi$ 
an embedding. Up to suspension, this image is thus the element
\[
[Z_n \to \R^n_+ \wedge_{\Diff} \Emb] \mapsto \MU^{n-2}(S^n B\Diff(F)) 
\]
defined by the tautological family $F \times_{\Diff} E\Diff(F)$ of surfaces 
over the classifying space of the diffeomorphism group. It is 
primitive in the Hopf-like structure defined by gluing: in fact it is
the image of 
\[
\sum_{k \geq 1} \kappa_k t_{k+1} \in \MU^{-2}(B\Gamma^+_\infty) \otimes \Q \;.
\]
If $v$ is a formal indeterminate of cohomological degree two, then the class
\[
\Phi = \exp (v\thom(-\eta)) \in \MU^0_\Q(\Omega^\infty \C P^\infty_{-1})[[v]]
\]
defined by finite unordered configurations of points on the universal surface 
(with $v$ a book-keeping indeterminate of cohomological degree two) is a kind
of exponential transformation
\[
\tilde \Phi_* : H_*(Q\C P^\infty_+,\Q) \to H_*(\MU,\Q[[v]]) \;.
\]
From this perspective it is natural to interpret the Thom class in 
$\MU^{-2}(\C P^\infty_{-1}) \otimes \Q$ as the sum $\sum_{k \geq -1} 
t_{k+1} \e^k$, with $t_0 = v^{-1}\e$. \bigskip

\noindent
{\bf 4.3} A class in the cohomology group
\[
H^2_{\rm Lie}(V,\R) \cong \Lambda^2(V^*)
\]
of a real vector space $V$ defines a Heisenberg extension
\[
0 \to \T \to H \to V \to 0 \;;
\]
the representation theory of such groups, and in particular the
construction of their Fock representations, is classical [5]. What is
important to us is that these are {\bf projective} representations
of $V$, with positive energy; such representations have very
special properties. \bigskip

\noindent
The loop group of a circle is a key example; it possesses an
intrinsic symplectic form, defined by formulae much like those of
\S 2 [23 \S 5, \S 7b]. Diffeomorphisms of the circle act on any such loop 
group, and it is a deep property of positive-energy representations, 
that they extend to representations of the resulting semidirect product
of the loop group by $\Diff S^1$. Therefore by restriction a positive-energy 
representation of a loop group automatically provides a representation 
of $\Diff S^1$. This [Segal-Sugawara [25 \S 13.4]] construction yields the
action of Witten's Virasoro algebra on the Fock space 
\[
{\rm Symm}(H^*(\C P^\infty_+)) \cong \Q[t_k \:|\; k \geq 0] \;.
\]
In Kontsevich's model, the classes $t_k$ are identified with the
symmetric functions
\[
\tr \; \gamma_{k-\half}(\Lambda^2) \sim - (2k-1)!! \; \tr \;
\Lambda^{-2k-1} 
\]
of a positive-definite Hermitian matrix $\Lambda$. \bigskip

\noindent
Note, however, that the deeper results of Kontsevich and Witten
theory [31] are inaccessible in this toy model: that theory is
formulated in terms of compactified moduli spaces $\M_g$ of
algebraic curves. The rational homology of $Q(\coprod \M_g)$
(suitably interpreted, for small $g$) contains a fundamental
class
\[
\exp(\sum_{g \geq 0} [\M_g]v^{3(g-1)})
\]
for the moduli space of not-necessarily-connected curves.
Witten's tau-function is the image of this `highest-weight'
vector under the analog of $\tilde \Phi$; it is killed by
the subalgebra of Virasoro generated by the operators $L_k$ 
with $k \geq -1$.  \bigskip

\section{Concluding remarks}

\noindent
{\bf 5.1} Witten has proposed a generalization of 2D topological
gravity which encompasses surfaces with higher spin structures:
for a closed smooth surface $F$ an $r$-spin structure is roughly
a complex line bundle $L$ together with a fixed isomorphism
$L^{\otimes r} \cong T_F$ of two-plane bundles, but for surfaces
with nodes or marked points the necessary technicalities are
formidable [14]. The group of automorphisms of such a structure is
an extension of its group of diffeomorphisms by the group of
$r$th roots of unity, and there is a natural analog of the group
completion of the category defined by such surfaces. The
generalized Madsen-Tillmann construction maps this loopspace to
the Thom spectrum $\thom(-\eta^r)$, and it is reasonable to expect
that this map is equivariant with respect to automorphisms of the
group of roots of unity. This fits with some classical homotopy 
theory: if (for simplicity) $r=p$ is prime, multiplication by an 
integer $u$ relatively prime to $p$ in the $H$-space structure of 
$\C P^\infty$ defines a morphism
\[
\thom(-\eta^p) \to \thom(-\eta^{up})
\]
of spectra, and the classification of fiber-homotopy equivalences
of vector bundles yields an equivalence of $\thom(-\eta^{up})$ with
$\thom(-\eta^p)$ after $p$-completion. There is an analogous
decomposition of $t_\T H\Z_p$ and a corresponding decomposition of
the associated Fock representations [20 \S 2.4]. \bigskip

\noindent
{\bf 5.2} Tillmann has also studied categories of surfaces above a parameter 
space $X$; the resulting group completions have interesting connections 
with both Tate and quantum cohomology. When $X$ is a compact smooth 
almost-complex manifold, its Hodge-deRham cohomology admits a natural
action of the Lie algebra generated by the Hodge dimension operator
$H$ together with multiplication by the first Chern class ($E$) and its 
adjoin ($F = *E*$) [26]. Recently Givental [9 \S 8.1] has shown that earlier 
work of (the schools of) Eguchi, Dubrovin, and others can be reformulated 
in terms of structures on $t^{*,*}_\T H_{\rm dg}(X)$, given a symplectic
structure generalizing that of \S 3. In this work, the relevant
involution is 
\[
I_{\rm Giv} = \exp(\half H) \exp(-E) \; I \; \exp(E) \exp(-\half H) \;;
\]
it would be very interesting if this involution could be understood 
in terms of the equivariant geometry of the free loopspace of $X$ [7].
\bigskip

\noindent
{\bf 5.3} Nothing forces us to restrict the construction of Madsen and
Tillmann to two-manifolds, and I want to close with a remark about the
cobordism category of smooth spin four-manifolds bounded by ordinary
three-spheres. A parametrized family of such objects defines, as in
\S 2.3, an element of the truncated equivariant cobordism group
\[
\tau^{-4}_{\SU(2)} \MSpin (X_+) \;.
\]
On the other hand, it is a basic fact of four-dimensional life that
\[
\Spin(4) = \SU(2) \times \SU(2) \;,
\]
so the Madsen-Tillmann spectrum for the cobordism category of such
spin four-folds is the twisted desuspension
\[
B\Spin(4)^{-\rho} = (\aitch P_\infty \times \aitch P_\infty)^{-V^* 
\otimes_{\aitch} V}
\]
of the classifying space of the spinor group by the representation $\rho$
defined by the tensor product of two standard rank one quaternionic 
modules over $\SU(2)$ [13 \S 1.4]. Composition with the Dirac operator 
defines an interesting rational homology isomorphism
\[
(\aitch P_\infty \times \aitch P_\infty)^{-V^* \otimes_{\aitch} V} \to
\aitch P_\infty^{-V} \wedge \MSpin \to \aitch P_\infty^{-V} \wedge \kO
\]
related in low dimensions to the classification of unimodular even indefinite 
lattices [27, 30]. This suggests that the Tate cohomology $t_{\SU(2)}^* \kO$
may have an interesting role to play in the study of topological gravity in
dimension four.  \bigskip

\newpage

\bibliographystyle{amsplain}

\begin{thebibliography}{99}

\bibitem [1]{1} J.F. Adams, \dots what we don't know about $\R P^{\infty}$, 
in {\bf New Developments in Topology}, ed. G. Segal,  LMS Lecture Notes 11 
(1972) 1 - 9

\bibitem [2]{2} A. Adem, R.L. Cohen, W. Dwyer,Generalized Tate homology, 
homotopy fixed points and the transfer, in {\bf Algebraic topology (Evanston
88)} 1 - 13, Contemp. Math. 96 (1989)

\bibitem [3]{3} V.I. Arnol'd, Symplectization, complexification and 
mathematical trinities, in {\bf The Arnoldfest} 23 - 37, Fields Inst. 
Commun. 24 (1999) 

\bibitem [4]{4} Th. Br\"ocker, E.C. Hook, Stable equivariant bordism, Math.
Zeits. 129 (1972) 269 - 277  

\bibitem [5]{5} P. Cartier, Quantum mechanical commutation relations and 
theta functions, in {\bf Algebraic Groups and Discontinuous Subgroups}, Proc. 
Sympos. Pure Math. 9 (1966) 361 - 383 

\bibitem [6]{6} R.L. Cohen, J.D.S. Jones, G.B. Segal, Floer's
infinite dimensional Morse theory and homotopy theory, in the {\bf Floer
Memorial Volume}, Birkh\"auser, Progress in Mathematics 133 (1995)
297 - 326

\bibitem [7]{7} --------, A. Stacey, Fourier decompositions of loop bundles,
Proc. Northwestern Conf (2001), to appear

\bibitem [8]{8} P. DiFrancesco, C. Itzykson, J.-B. Zuber, Polynomial 
averages in the Kontsevich model, CMP 151 (1993) 193-219

\bibitem [9]{9} A. Givental, Gromov - Witten invariants and quantization
of quadratic hamiltonians, available at {\tt math.AG/0108100}

\bibitem [10]{10} J.P.C. Greenlees, A rational splitting theorem for the 
universal space for almost free actions, Bull. London Math. Soc. 28 (1996) 
183 - 189

\bibitem [11]{11} --------, J.P. May, {\bf Generalized Tate cohomology}, 
Mem. AMS 113 (1995) 

\bibitem [12]{12} --------, H. Sadofsky, The Tate spectrum of $v_n$-periodic
complex-oriented theories, Math. Zeits. 222 (1996) 391 - 405

\bibitem [13]{13} R. Gompf, A. Stipsicz {\bf Four-manifolds and Kirby 
calculus}, AMS Grad Texts 20 (1999)

\bibitem [14]{14} T. Jarvis, T. Kimura, A. Vaintrob, Moduli spaces
of higher spin curves and integrable hierarchies, available at
{\tt math.AG/9905034}

\bibitem [15]{15} T. J\'ozefiak, Symmetric functions in the Kontsevich-Witten 
intersection theory of the moduli space of curves, Lett. Math. Phys. 33
(1995) 347-351 

\bibitem [16]{16} M. Kapranov, E. Vasserot, Vertex algebras and the formal 
loop space, available at {\tt math.AG/0107143}

\bibitem [17]{17} I. Madsen, U. Tillmann,  The stable mapping-class group
and $Q(\C P^{\infty}_+)$, Invent. Math. 145 (2001) 509 - 544.

\bibitem [18]{18} M. Mahowald, {\bf On the metastable homotopy of} $S^{n}$, 
Mem. AMS 72 (1967)

\bibitem [19]{19} J. Morava, Cobordism of involutions revisited, revisited, in
{\bf The Boardman Festschrift}, Contemporary Math. 239 (1999)

\bibitem [20]{20} --------, An algebraic analog of the Virasoro group,
Czech. J. Phys. 51 (2001), available at {\tt math.QA/0109084}

\bibitem [21]{21} D. Quillen, On the formal group laws of unoriented
and complex cobordism theory, BAMS 75 (1969) 1293 - 1298

\bibitem [22]{22} --------, Elementary proofs of some results of cobordism
theory using Steenrod operations, Adv. in Math 7 (1971) 29 - 56

\bibitem [23]{23} G. Segal, Unitary representations of some 
infinite-dimensional groups, Comm. Math. Phys. 80 (1981) 301 - 342. 

\bibitem [24]{24} -------, A. Pressley, {\bf Loop groups}, Oxford (1986)

\bibitem [25]{25} -------, {\bf Algebres de Lie semisimples complexes},
Benjamin (1966)

\bibitem [26]{26} J.P. Serre, {\bf Corps Locaux}, Hermann (1968)

\bibitem [27]{27} -------, {\bf A Course in Arithmetic}, Springer (1973)

\bibitem [28]{28} R. Swan, Periodic resolutions for finite groups, 
Annals of Math. 72 (1960) 267 - 291

\bibitem [29]{29} J. Tate, Residues of differentials on curves, Ann. Sci. 
Ecole Norm. Sup. 1 (1968) 149 - 159

\bibitem [30]{30} C.T.C. Wall, On simply-connected four-manifolds, J.
London Math. Soc. 39 (1964) 141 - 149 

\bibitem [31]{31} E. Witten, Two-dimensional gravity and intersection theory 
on moduli space, Surveys in Differential Geometry 1 (1991) 243 - 310

 
\end{thebibliography}

\end{document}